\documentclass{psapm-l}

\newtheorem{theorem}{Theorem}[section]

\theoremstyle{definition}

\theoremstyle{remark}
\newtheorem{remark}[theorem]{Remark}

\numberwithin{equation}{section}

\newcommand{\ee}{\varepsilon}
\newcommand{\ue}{u^{\varepsilon}}
\newcommand{\zz}{\zeta}

\begin{document}
\bibliographystyle{plain}
\title[Invariant manifold for viscous profiles]{Invariant Manifolds for Viscous Profiles of a Class of Mixed Hyperbolic-Parabolic Systems}

\author{Stefano Bianchini}
\address{SISSA, via Beirut 2-4 34014 Trieste, Italy}
\email{bianchin@sissa.it}

\author{Laura V. Spinolo}
\address{Centro di Ricerca Matematica Ennio De Giorgi, Scuola Normale Superiore, 
Piazza dei Cavalieri 3, 
56126 Pisa, 
Italy}
\email{laura.spinolo@sns.it}


\begin{abstract}
We are concerned with viscous profiles (travelling waves and steady solutions) for mixed hyperbolic-parabolic systems in one space variable. For a class of systems including the compressible Navier Stokes equation, these profiles satisfy a singular ordinary differential equation in the form
\begin{equation}
\label{e:ab}
     \frac{d U}{d t} = \frac{1}{\zeta (U)} F(U).
\end{equation}
Here $U$ takes values in $\mathbb R^d$ and $F: \mathbb R^d \to \mathbb R^d$ is a regular function. The real valued function $\zeta(U)$ is as well regular, but the equation is \emph{singular} because $\zeta (U)$ can attain the value $0$. We focus on a small enough neighbourhood of a point $\bar U$ satisfying $F(\bar U) = \vec 0$, $\zeta (\bar U) =0$. From the point of view of the applications to the study of hyperbolic-parabolic systems this means restricting to systems with small total variation.

We discuss how to extend the notions of center manifold and of uniformly stable manifold. Also, we give conditions ensuring that if $\zeta (U) \neq 0$ at $t=0$ then $\zeta (U) \neq 0$ at every $t$. We provide an example showing that if $\zeta(U)$ becomes  zero in finite time then in general the solution $U$ of equation \eqref{e:ab}  is not continuously differentiable.
\end{abstract}

\maketitle

\section{Viscous profiles for mixed hyperbolic-parabolic systems}
We are concerned with systems in the form 
\begin{equation}
\label{e:p}
        E(u) u_t + A(u, \, u_x) u_x = B (u) u_{xx}.
\end{equation}
Here the function $u$ takes values in $\mathbb{R}^N$ and depends on two scalar variables, $t$ and $x$. The matrices $E$, $A$ and $B$ have all dimension $N \times N$. 
The conservative case 
\begin{equation*}
          u_t + f(u)_x  = \Big( B(u) u_x \Big)_x  
\end{equation*}
is, in particular, included in the formulation \eqref{e:p}. 

In the following, we will focus on mixed hyperbolic-parabolic systems, i.e. we will assume that matrix $B$ in \eqref{e:p} is singular. This case case is interesting from the point of view of applications. Indeed, it is the case of the compressible Navier Stokes equation in one space variable: 
\begin{equation}
\label{e:ns:eul}
       \left\{
       \begin{array}{lll}
              \rho_t + ( \rho v )_x =0 \\
	      (\rho v)_t + \Big( \rho v^2 + p \Big)_x = \displaystyle{ \Big( \nu v_x  \Big)_x } \\
	      \displaystyle{ \Big( \rho e + \rho \frac{v^2}{2}\Big)_t + \Big(v \Big[ \frac{1}{2} \rho v^2 
	      + \rho e + p \Big] \Big)_x = \Big( k \theta_x + 
	      \nu v v_x \Big)_x}. \\
       \end{array}
       \right.
\end{equation}
Here the unknowns are $\rho(t, \, x), \, v(t, \, x)$ and $e(t, \, x)$: the function $\rho$ is the density of the fluid, $v$ represents the velocity of the particles in the fluid and $e$ is the internal energy. The function $p= p(\rho, \, e) >0$ is the pressure and satisfies $p_{\rho} >0$, while $\theta$ is the absolute temperature and in the case of a polytropic gas satisfies
$$
   \theta = \frac{e (\gamma -1 )}{R},  
$$     
where $R$ is the universal gas constant and $\gamma$ is a constant specific of the gas. Finally, $\nu(\rho)>0$ and $k(\rho)>0$ are the viscosity and the heat conduction coefficients respectively. 

In \cite{KawShi:normal}, Kawashima and Shizuta introduced a set of hypotheses that are satisfied by the equations of the hydrodynamics and of the magnetohydrodynamics and that are frequently exploited to study the hyperbolic-parabolic system \eqref{e:p}. In the following we will assume that the Kawashima Shizuta conditions are satisfied but, since we do not exploit them explicitly, we do not recall them. 

We are concerned with special classes of solutions of \eqref{e:p}, namely travelling waves and steady solutions. Travelling waves satisfy the ordinary differential equation 
\begin{equation}
\label{e:tw}
          \big[ A ( U, \, U') - \sigma E (U) \big] U ' = B (U) U '',
\end{equation}          
while steady solutions are solutions of the ODE 
\begin{equation}          
\label{e:bl}
            A ( U, \, U')  U ' = B (U) U '' . 
\end{equation}
In \eqref{e:tw}, $\sigma$ is a real constant which is usually called the speed of the wave. Note that from a solution $U(y)$ of \eqref{e:tw} we can obtain a solution of the original hyperbolic-parabolic system \eqref{e:p} by setting $u(t, \, x) = U(x - \sigma t)$. Moreover, any solution $U(x)$ of \eqref{e:bl} is a steady solution of \eqref{e:p}, i.e. a solution that does not depend on time. Also, in the following we will focus on the case $\sigma$, the speed of the travelling wave \eqref{e:tw}, is close to an eigenvalue of the matrix $A(U, \, \vec 0)$. Since in general $0$ is not an eigenvalue of $A(U, \, \vec 0)$, we keep the cases \eqref{e:tw} and \eqref{e:bl} separated. 

It is known that the study of travelling waves and steady solutions can provide useful information to study the limit 
$\ee \to 0^+$ of the family of functions $\ue$ satisfying 
\begin{equation}
\label{e:visc}
     E(\ue) \ue_t + A(\ue, \, \ee \ue_x) \ue_x = \ee B (\ue) \ue_{xx}. 
\end{equation}
The literature concerning this issue is very wide, so we just refer  to Benzoni-Gavage, Rousset, Serre and Zumbrun \cite{BenRouSerreZum},  to Rousset \cite{Rousset:char}, to Zumbrun \cite{Zum:review}, and to the rich bibliography contained therein. For a more general introduction to the parabolic approximation of hyperbolic problems we refer instead to the books by Dafermos \cite{Daf:book} and by Serre \cite{Serre:book} and to the references therein. 

Note that, when system \eqref{e:p} is the Navier Stokes equation, the system we obtain formally setting $\ee =0$ in \eqref{e:visc} is the Euler equation.  

In studying the viscous profiles of the Navier Stokes equation \eqref{e:ns:eul} one encounters a singular ordinary differential equation in the form 
\begin{equation}
\label{e:sin}
     \frac{dU}{dt } = \frac{1}{\zz (U)} F(U).
\end{equation}
Here the unknown $U$ is vector valued and has the same dimension as the function $F$. The function $\zz$ is real valued and the singularity of the equation 
comes from the fact that $\zeta$ can attain the value $0$. The link between \eqref{e:ns:eul} and \eqref{e:sin} is the following. Assume that we want to study the steady solutions of the Navier Stokes equation: we thus focus on \eqref{e:bl}.  Set  
\begin{equation}
\label{e:h:ns:ux}
   w = \rho_x  \qquad \vec z = \Big( v_x, \, e_x  \Big)^t
\end{equation}
After some computations (see Bianchini and Spinolo \cite{BiaSpi:ode} for the details), we get that the steady solutions of the Navier Stokes equation satisfy 
\begin{equation}
\label{e:matrix}
      \left(
             \begin{array}{cc}
                   a_{11} v  & A_{21}^t \\
                   A_{21} & A_{22} \\
             \end{array}
             \right) 
             \left(
             \begin{array}{ll}
                         w \\
                         \vec z \\
             \end{array}
             \right)
             = 
             \left(
             \begin{array}{cc}
                   0 & 0 \\
                   0 & b(u) \\
             \end{array}
             \right) 
             \left(
             \begin{array}{ll}
                         w_x \\
                         \vec z_x \\
             \end{array}
             \right).   
\end{equation}
Here $A_{22}$ and $b$ are $2 \times 2$ matrices, $a_{11}$ is a real valued function, the function $A_{21}$ takes values into $\mathbb R^2$ and  $A_{21}^t$ denotes its transpose. The exact expression of these terms is not important here: we just point out that the matrix $b$ is invertible. 
Equation \eqref{e:matrix} gives 
$$
  \left\{
\begin{array}{ll}
           a_{11}  v  w                 +    A_{21}^t \vec  z = 0 \\
            A_{21} w        +    A_{22} \vec z = b \vec z_x
\end{array}           
\right.
$$
Assume $a_{11} v \neq 0$, then we get
\begin{equation}
\label{e:h:ns:zx}
\left\{
\begin{array}{ll}
            w =  \displaystyle{    -  \frac{ A_{21}^t \vec{z} }{a_{11}   v}     } \\
            \vec z_x =   b^{-1} \displaystyle{\Big[     A_{22}   -   \frac{A_{21} A_{21}^t   }{ a_{11} v}   \Big] \vec z } \\
\end{array}
\right.
\end{equation}
Note that the previous expression is well defined since the matrix $b$ is invertible. Summing up, we get that the steady solutions of the Navier Stokes equation can be written in the form \eqref{e:sin} provided that $U = ( \rho, \, v, \, e, \, \vec z )^t$, $\zz (U) = v$ and 
$$
    F(U) =
    \left(
    \begin{array}{ccc}
                A_{21}^t \vec{z} / a_{11} \\
                v \, \vec z   \\
                 b^{-1} \displaystyle{\Big[     A_{22} v   -   A_{21} A_{21}^t  /a_{11}   \Big] \vec z } \\
    \end{array}
    \right).
$$
Some remarks are here in order: first, the exact expression of the function $a_{11}$ is $p_{\rho}/ \rho^2$, where $p_{\rho}>0$ is the partial derivative of the pressure with respect to the variable $\rho$. We can restrict to the case $\rho$ is strictly positive and bounded away from zero: this implies that vacuum states are not assumed. We then have that the function $a_{11}$ is well defined and does not attain the value zero. On the other side, $v$ represents the velocity of the fluid and in general it can attain the value zero, which is the singular value for the equation satisfied by the steady solutions. 
 
 Moreover, we underline that so far we have considered only steady solutions. However, also the equation of travelling wave profiles \eqref{e:tw} may become singular, in the sense of \eqref{e:sin}. 
 
 Also, the considerations carried on so far can be extended to a larger class of mixed hyperbolic-parabolic systems that do not satisfy a condition of block linear degeneracy defined in Bianchini and Spinolo \cite{BiaSpi:rie}.  
 
 \begin{remark}
 The reason why we have the term $a_{11}$ in~\eqref{e:matrix} is the following: instead of working directly on~\eqref{e:ns:eul}, we consider an equivalent symmetric system in the form~\eqref{e:p} (see again Bianchini and Spinolo \cite{BiaSpi:ode} for the details). System \eqref{e:p} is \emph{symmetric} if  
 $$
      A(u, \, \vec 0)^t = A(u, \vec 0). 
 $$
 \end{remark}
\section{Invariant manifolds for a class of singular ODEs}
In the following we will thus focus on the singular ordinary differential equation
\begin{equation}
\label{e:sin2}
    \frac{d U }{d t} = \frac{1}{ \zz (U)} F(U). 
\end{equation}
We will be concerned with the solutions $U$ belonging to a small enough neighbourhood of a value $\bar U$ satisfying $\zz (\bar U) =0$, $F(\bar U) =\vec 0$. From the point of view of the applications to the analysis of system \eqref{e:p}, this implies that we are restricting to steady solutions and travelling waves having small enough total variation. Without loss of generality, in the following we will assume that $\bar U= \vec 0$. 

What we are interested in is the existence of locally invariant manifolds for \eqref{e:sin2}. If $U$ takes values in $\mathbb R^d$, a locally invariant manifold $\mathcal M$ for \eqref{e:sin2} is contained in $\mathbb R^d$ and satisfies the following: if $u_0 \in \mathcal M$, then the solution of the Cauchy problem
\begin{equation*}
\left\{
\begin{array}{ll}
           \displaystyle{ \frac{d U }{d t} = \frac{1}{ \zz (U)} F(U)} \\
           u (0) = u_0 \\
\end{array}
\right.
\end{equation*}
belongs to $\mathcal M$ is $|t|$ is small. In particular, we are interested in extending the notions of \emph{uniformly stable} and of \emph{center manifold} to the case of the singular ODE \eqref{e:sin2}. We recall here that a center manifold for the non singular ODE
\begin{equation}
\label{e:nonsin}
            \frac{d U }{d t} =  G(U) \qquad U \in \mathbb R^d 
\end{equation}
is defined in a neighbourhood of an equilibrium point $\bar U$. Loosely speaking,  a \emph{center manifold} contains the orbits of \eqref{e:nonsin} that are globally bounded for $t \to \pm \infty$ , more precisely $|U(t)| \leq \delta $ for every $t$. The constant $\delta$ is the size of the neighbourhood and depends on the function $G$ and on the point $\bar U$. We refer to the book by Katok and Hasselblatt \cite{KHass} for a complete discussion. Also, a presentation of the the most important properties of the notion of center manifold is given in Bressan's notes \cite{Bre:note}. For the applications of the notion of center manifold to the study of the parabolic approximation of hyperbolic problems see for example Bianchini and Bressan \cite{BiaBrevv}. 

The \emph{stable manifold} of \eqref{e:nonsin} contains the orbits that when $t \to + \infty$ converge exponentially fast to the equilibrium point $ \bar U$. We refer to the book by Perko \cite{Perko} for a complete discussion.  The notion of \emph{uniformly stable manifold} is an extension of the notion of stable manifold. Assume that $E$ is a manifold of equilibria for system \eqref{e:nonsin} and assume that $\bar U \in E$. The uniformly stable manifold relative to $ E$ contains the orbits of system \eqref{e:nonsin} that  for $t \to + \infty$ decay eponentially fast to an equilibrium point in $E$. It is defined in a small enough neighbourhood of $\bar U$ and contains the stable manifold. However, in general the inclusion is strict since the orbits on the stable manifold can converge only to $\bar U$, while when we are on the uniformly stable manifold the limit can vary on $E$. The uniformly stable manifold is sometimes called the \emph{slaving manifold} relative to $E$. Its existence can be viewed as a consequence of Hadamard-Perron Theorem: we refer again to the book by Katok and Hasselblatt \cite{KHass} for the Hadamard-Perron Theorem, while a specific discussion on the uniformly stable manifold and its applications to the study of the parabolic approximation of  hyperbolic systems can be found in Ancona and Bianchini \cite{AnBia}. 

It is known that, from the point of view of the applications to the study of the parabolic approximation \eqref{e:visc}, if we restrict to systems having small enough total variation it is interesting to focus on travelling waves lying on a center manifold and steady solutions lying on a stable or on a uniformly stable manifold.  

There is a rich literature concerning the family of systems 
\begin{equation}
\label{e:i:spt}
          \frac{d U^{\ee} }{d t} = \frac{1}{ \ee } \,  F(U, \, \ee) \qquad U \in \mathbb R^d.
\end{equation}
Here we just refer to the notes by Jones \cite{Jones} and to the rich bibliography contained therein. In particular, \cite{Jones} provides a nice overview of some papers by Fenichel (\cite{Fenichel, Fenichel_per} for example),  whose ideas are exploited in the following.  In \eqref{e:i:spt} the singularity $\ee$ is a parameter, $\ee \to 0^+$. The main novelty here is that in~\eqref{e:sin2} we consider the case $\zz(U)$ is a nontrivial function of the solution $U$ itself. As a consequence, we have to take into account the possibility that $\zz (U) \neq 0$ when $t =0$, but $\zz (U)$ attains the value $0$ in finite time. This may lead to a loss of regularity in the solution $U$. Consider, for example, the following system: 
  \begin{equation}
           \label{e:ex:fast}
           \left\{
           \begin{array}{ll}
                      d u_1 / dt = - u_2 / u_1 \\
                      d u_2 / dt = - u_2 . \\
           \end{array}
           \right.                     
 \end{equation}    
It can be written in form \eqref{e:sin2} provided that $U= (u_1, \, u_2)^t$, $\zeta (U) = u_1$ and 
$$
    F (U) =
    \left(
    \begin{array}{cc}
                - u_2 \\
                -u_2 u_1  \\
     \end{array}
     \right)
$$
The solution of \eqref{e:ex:fast} is 
\begin{equation}
\label{e:ex:sol}
     \left\{
    \begin{array}{lll}
                \displaystyle{ u_1 (t) = \sqrt{  u_1 (0) +  u_2 (0) \big(  e^{- t } - 1 \big)  }} \\
                \\
                   \displaystyle{  u_2 (t) = u_2 (0) e^{- t } } \\
    \end{array}
    \right.
\end{equation}
Choosing $u_2 (0) > u_1 (0) >0$, one has that $\zz (U) =u_1 (t)$ can attain the singular value $0$ for a finite $t$. Note that at that point $t$ the first derivative $d u_1 / dt$
blows up: thus, in particular, the solution \eqref{e:ex:sol} of \eqref{e:ex:fast} is not $\mathcal C^1$. 

In the following, we will look for conditions that rule such a loss of regularity: this sounds reasonable in view of the applications to the analysis of the viscous profiles. Indeed, when considering the parabolic approximation \eqref{e:visc} it seems reasonable to look for regular solutions. To prevent losses of regularity like the one in \eqref{e:ex:sol} we need to be sure that if a solution $U$ of \eqref{e:sin2} satisfies $\zz (U) \neq 0$ at $t=0$, then $\zz (U) \neq 0$ for every $t$. In the following, we will state conditions ensuring that this property holds. 
\begin{remark}
           To simplify the notations, in the following we will always assume that $\zz (U)>0$ at $t=0$. The case $\zz (U) <0$ does not involve additional difficulties and can be tackled with techniques similar to those discussed here. 
\end{remark}
\section{Main results}
Summing up, our goals are as follows: we want to define locally invariant manifolds that extend the definition of center and uniformly stable manifold and we want to make sure that losses of regularity like the one in \eqref{e:ex:sol} are ruled out.

We state our main results here, postponing to Section \ref{s:hyp} the statement and the discussion of the hypotheses. We refer to \cite{BiaSpi:ode} for a proof of Theorems~\ref{t:center} and~\ref{t:stable}. Concerning the definition of a center manifold, we have the following.  
\begin{theorem}
\label{t:center}
            Let all the hypotheses introduced in Section \ref{s:hyp} hold. Also, assume that 
            $$
                 F( \vec 0) = \vec 0 \qquad \zz (\vec 0) = 0.
            $$
            Then there exists a manifold $\mathcal M^c$, defined in a small enough neighbourhood of $\vec 0$ and satisfying the following properties. The manifold $\mathcal M^c$ is locally invariant for 
            \begin{equation}
            \label{e:sin3}
                 \frac{d U }{d t} = \frac{1}{ \zz (U)} F(U).
            \end{equation}
            It contains all the solutions $U$ satisfying $|U(t)| \leq \delta$ for every $t$. Here $\delta >0$ is a constant depending on the system. Also, if $U(t)$ is an orbit lying on $\mathcal M^c$ and $\zz (U) \neq 0$ at $t=0$, then $\zz (U) \neq 0$ for every $t$. 
\end{theorem}
Before discussing the extension of the notion of uniformly stable manifold, some remarks are in order. One can show (see Bianchini and Spinolo \cite{BiaSpi:ode} for the details) that the solutions of \eqref{e:sin3} that lie on the manifold $\mathcal M^c$ satisfy an equation which actually has no singularity in it. The situation is completely different if we consider the uniformly stable manifold. To see this, let us consider a trivial example:
  \begin{equation}
           \label{e:ex:triv}
           \left\{
           \begin{array}{ll}
                      d u_1 / dt = -  5 u_1  \\
                      d u_2 / dt = - u_2 / \ee \\
                      d \ee / dt = 0. 
           \end{array}
           \right.                     
 \end{equation}   
Note that the subspace $E = \{ (0, \, 0, \, \ee ): \; \ee \in \mathbb R \}$ is entirely made by equilibria of \eqref{e:ex:triv}. If $\ee >0$ the uniformly stable manifold relative to $E$ is the whole space $\mathbb R^3$, since any solution of \eqref{e:ex:triv} decays exponentially fast to a point in $E$. However, the first component goes like $e^{- 5 t}$: thus, the speed of exponential decay is bounded in $\ee$ and $u_1$ is not affected by the presence of the singularity. Conversely, the second component decays like $e^{-t / \ee}$ and hence the speed of exponential decay gets faster and faster as $\ee \to 0^+$. The second component of the solution is thus strongly affected by the presence of the singularity and can be regarded as a \emph{fast dynamic}, while the first component is a \emph{slow dynamic}.  Summing up, in \eqref{e:ex:triv} any orbit lying on the uniformly stable manifold relative to $E$ decomposes as the sum of a slow and a fast dynamic. 

This behaviour is somehow inherited by general non linear systems, in the following sense. If we perform the change of variable $\tau = t / \ee$, system \eqref{e:ex:triv} becomes  
$$
      \left\{
           \begin{array}{ll}
                      d u_1 / d \tau = -  5 u_1 \ee  \\
                      d u_2 / d \tau = - u_2  \\
                      d \ee / d \tau = 0. 
           \end{array}
           \right.     
$$
We can thus single out the \emph{fast dynamics} by saying that they are solutions of \eqref{e:ex:triv} that are exponentially decaying with respect to the variable $\tau$. In the general non linear case, we consider system 
\begin{equation}
\label{e:nsin}
     \frac{d U}{ d \tau} = F(U), 
\end{equation}
which is \emph{formally} obtained from system \eqref{e:sin3} through the change of variables ${\tau = \tau (t)}$ defined by
\begin{equation}
\label{e:i:changev}
\left\{
\begin{array}{lll}
          \displaystyle{   \frac{d \tau}{dt} = \frac{1}{\zeta [U(t) ]} } \\
         \\
          \tau ( 0) =0.\\
\end{array}
\right.
\end{equation}
The \emph{fast dynamics} of \eqref{e:sin3} are then the solutions of \eqref{e:nsin} that are exponentially decaying to zero with respect to the variable $\tau$. Note that, a priori, the change of variable \eqref{e:i:changev} is not well defined, because it may happen that $\zz (U)$ attains the value $0$ for a finite value of $t$ and hence that the solution of the Cauchy problem \eqref{e:i:changev} is not a diffeomorphism ${\tau: [0, \, + \infty[ \to [0, \, + \infty[ }$. 

Conversely, the \emph{slow dynamics} behave somehow like the component $u_1$ in \eqref{e:ex:triv}, namely they are solutions of \eqref{e:sin3} that satisfy an equation with no singularity in it. We refer to Bianchini and Spinolo \cite{BiaSpi:ode} for a more rigourous definition of fast and slow dynamics.

Before introducing Theorem \ref{t:stable}, we point out that, as a consequence of the hypotheses introduced in Section \ref{s:hyp}, we can perform a local change of coordinates such that, in the new coordinates, $\zz (U)$ is actually one of the components of $U$, say $\zz (U) = u_1$ if $U=(u_1, \dots u_d )^t.$ Also, the manifold ${E = \{ (u_1, \, 0 \dots 0): \; u_2 = u_3 = \dots u_d =0 \}}$ is entirely made by equilibria for 
 \begin{equation}
            \label{e:sin4}
                 \frac{d U }{d t} = \frac{1}{ \zz (U)} F(U).
            \end{equation}
Our result is the following:
\begin{theorem}
\label{t:stable}
          Assume that 
            $$
                 F( \vec 0) = \vec 0 \qquad \zz (\vec 0) = 0
            $$
            and that all the hypotheses introduced in Section \ref{s:hyp} are satisfied.
            Then there exists a manifold $\mathcal M^{us}$, defined in a small enough neighbourhood of $\vec 0$ and satisfying the following properties. 
            \begin{enumerate}
            \item             
            The manifold $\mathcal M^{us}$ is locally invariant for \eqref{e:sin4} and 
            it contains all the solutions $U(t)$ decaying exponentially fast to an equilibrium point in \\
                        {$E = \{ (u_1, \, 0 \dots 0): \; u_2 = u_3 = \dots u_d =0 \}$} . 
            \item
            Any orbit lying on 
            $\mathcal   M^{us}$ decomposes as 
            $$
                 U(t) = U^{slow}(t) + U^{fast} (t) + U^{pert} (t).  
            $$
            Here $ U^{slow}(t)$ and $ U^{fast}(t)$ are a slow and a fast dynamic respectively, in the sense explained before. The perturbation term $ U^{pert} (t)$ is due to the non linearity and satisfies 
            $$
                 | U^{pert} (t)| \leq C \; U^{fast} (0) \,  \zz \Big (U (0) \Big) \, 
             $$
            for a suitable constant $C$. 
            \item If $U(t)$ is an orbit lying on $\mathcal M^c$ and $\zz (U) > 0 $ at $t=0$, then $\zz (U) > 0$ for every $t$. Also, the Cauchy problem \eqref{e:i:changev} defines a diffeomorphism ${\tau: [0, \, + \infty [ \to  [0, \, + \infty [ .}$
            \end{enumerate} 
\end{theorem} 
\begin{remark}
By direct check one can verify that all the hypotheses introduced in Section~\ref{s:hyp} are verified by~\eqref{e:h:ns:zx}, the equation satisfied by the steady solutions of the Navier Stokes written in Eulerian coordinates. 

The application of Theorems~\ref{t:center} and ~\ref{t:stable} to~\eqref{e:h:ns:zx} is connected to a remark due to Fr\'{e}deric Rousset, which is the following. In our hypotheses the Lagrangian and the Eulerian coordinates are equivalent. However, steady solutions of the Navier Stokes equation written using Lagrangian coordinates are regular since they satisfy an ODE with no singularity. As shown by example~\eqref{e:ex:triv}, in general a solution of a singular ODE is not $\mathcal C^1$. If the solutions of~\eqref{e:h:ns:zx} were not $\mathcal C^1$, this would contradict the equivalence between Eulerian and Lagrangian coordinates. 

However, Theorems~\ref{t:center} and~\ref{t:stable} tell us that, if we restrict to solutions lying on either $\mathcal M^c$ or $\mathcal M^{us},$ then any loss of regularity is ruled out: if $v(0) >0$, then $v(t) >0$ for every $t$.   

For a different approach to the analysis of the viscous profiles of the Navier Stokes equation in Eulerian coordinates see for example Wagner \cite{Wagner} and the references therein. 
\end{remark}
\section{Hypotheses}
\label{s:hyp}
In this section, we introduce the hypotheses we impose on system \eqref{e:sin4}. All these conditions are satisfied by \eqref{e:h:ns:zx}, the equation satisfied by the viscous profiles  of the compressible Navier Stokes equation in one space variable written in Eulerian coordinates.

Concerning the regularity of the maps $F$ and $\zz$, we assume that they are $\mathcal C^{3}$. Also, without loss of regularity we can restrict to the case the equilibrium $\bar U$ is $\vec 0$, namely
$$
          F( \vec 0) = \vec 0 \qquad \zz (\vec 0) = 0
$$
We also assume the following conditions:
\begin{enumerate}
\item
         The gradient $\nabla \zeta (\vec 0)  \neq \vec 0$. 
\end{enumerate}
Let $\mathcal S$ be the singular set
\begin{equation*}
         \mathcal S : = \big\{ U: \; \zeta (U) =0   \big\}.
\end{equation*}
Thanks to the implicit function theorem, Hypothesis~(1) ensures that in a small enough neighbourhood of $\vec 0$ the set $\mathcal S$ is actually a manifold of dimension $d-1$, where $d$ is the dimension of $U$.  
\begin{enumerate} \setcounter{enumi}{1}
\item
          Let $\mathcal M^c$ be any center manifold for 
          \begin{equation}
           \label{e:nsin2}
                    \frac{d U}{ d \tau} = F(U)
          \end{equation}
          around the equilibrium point $\vec 0$. If $|U|$ is sufficiently small and $U$ belongs to the intersection $\mathcal M^c \cap \mathcal S$ , then $U$ is an equilibrium for  \eqref{e:nsin2}, namely $F(U)= \vec 0$ . 
\end{enumerate}          
 The reason why we introduce Hypothesis~(2) is the following. Consider the linear system 
 \begin{equation}
 \label{e:giro}
       \left\{
           \begin{array}{ll}
                      d u_1 / d \tau =  u_2 / \ee  \\
                      d u_2 / d \tau = - u_1 / \ee   \\
                      d \ee / d \tau = 0. 
           \end{array}
           \right.   
 \end{equation}
 The first component of the solution is 
 $$
     u_1 (t) = A \cos (t / \ee) + B \sin (t / \ee), 
 $$
 where $A$ and $B$ are real parameters. Letting $\ee \to 0^+$, we get that in general if $t \neq 0$ there exists no pointwise limit of $u_1$. Note that~\eqref{e:giro} does not satisfy Hypothesis~(2). Indeed, the center space is the whole $\mathbb R^3$. However, 
 $$
     F(U) =   
     \left(
           \begin{array}{ccc}
                        u_2  \\
                     - u_1    \\
                     0
           \end{array}
           \right)   
 $$ 
 is not identically zero when $\ee =0$.          
\begin{enumerate} \setcounter{enumi}{2}            
\item 
         There exists  a manifold of equilibria  $\mathcal M^{eq}$ for \eqref{e:nsin2} which contains $\vec 0$ and which is transversal to $\mathcal S$.
\end{enumerate}
Let $n_{eq}$ be the dimension of $\mathcal M^{eq}$. We recall that the manifolds $\mathcal S$ and $\mathcal M^{eq}$ are transversal if the intersection $\mathcal S \cap \mathcal M^{eq}$ is a manifold with dimension $n_{eq} - 1$ (as pointed out before, the dimension of $\mathcal S$ is $d-1$).   
\begin{enumerate} \setcounter{enumi}{3}
\item
          For every $U \in \mathcal S$, 
          \begin{equation*}
                    \nabla \zeta (U) \cdot F(U) =0  .   
          \end{equation*}
\end{enumerate}
Hypothesis~(4) is necessary if we want to rule out losses of regularity like the one in \eqref{e:ex:fast}. Indeed, example \eqref{e:ex:fast} satisfies all the hypotheses stated in this section, but Hypothesis~(4).

Thanks to Hypothesis~(4) and to the regularity of the functions $\zeta$ and $F$, the function
$$
    G(U) = \frac{   \nabla \zeta (U) \cdot F(U)}{ \zeta (U)}
$$
can be extended and defined by continuity on the surface $\mathcal S$. 
\begin{enumerate} \setcounter{enumi}{4}
\item
          Let $U \in \mathcal S$ be an equilibrium for \eqref{e:nsin2}, namely $\zeta (U)= 0$ and $F(U) = \vec 0$. Then 
          \begin{equation*}
                            G(U) = 0. 
          \end{equation*}
\end{enumerate}

\bibliography{biblio}
\end{document}